\documentclass[12pt]{article}
\usepackage{amssymb,amsmath,amsthm}

\newtheorem{thm}{Theorem}
\newtheorem{lem}{Lemma}
\newtheorem{defn}{Definition}

\newtheorem{rem}{Remark}
\newtheorem{cor}{Corollary}

\newcommand {\C} {{\mathbb C}}
\newcommand {\N} {{\mathbb N}}
\newcommand {\Z} {{\mathbb Z}}
\newcommand {\R} {{\mathbb R}}

\newcommand {\T} {{\mathbb T}^{2}}
\newcommand {\Td} {{\mathbb T}^{d}}

\begin{document}
\title{Higher Order Oscillating Sequences, Affine Distal Flows on the $d$-Torus, and Sarnak's Conjecture~~}

\date{}
\author{Yunping JIANG\footnote{The author is partially supported by the collaboration grant from the Simons Foundation [grant number 199837] and the CUNY collaborative incentive research grants [grant number 2013] and awards from PSC-CUNY and grants from NSFC [grant numbers 11171121 and 11571122].}}

\maketitle

\abstract{In this paper, we give two precise definitions of a higher order oscillating sequence 
and show the importance of this concept in the study of Sarnak's conjecture. 
We prove that any higher order oscillating sequence of order $d$ is linearly disjoint 
from all affine distal flows on the $d$-torus for all $d\geq 2$. One consequence of this result is that 
any higher order oscillating sequence of order $2$ is linearly disjoint 
from all affine flows on the $2$-torus with zero topological entropy. 
In particular, this reconfirms Sarnak's conjecture for all affine flows 
on the $2$-torus with zero topological entropy and 
for all affine distal flows on the $d$-torus for all $d\geq 2$. }

\section{Introduction}
Suppose $X$ is a compact metric space with metric $d (\cdot, \cdot)$. Let $T: X\to X$ be a continuous map. We call $T$ a {\em flow} or a {\em dynamical system} because we will consider iterations $\{ T^{n}\}_{n=0}^{\infty}$. Let $\C$ denote the complex plane. Let $\R$,
$\R^{2}$, and $\R^{d}$ denote the real line, the real plane, and the $d$-Euclidean space, respectively. We denote by $\N$ the set of positive integers and by $\Z$ the set of integers. Then $\Z^{2}$ and $\Z^{d}$ are the integer lattices in $\R^{2}$ and $\R^{d}$. Let $C(X, \C)$ be the space of all continuous functions $f: X\to \C$. 

Suppose ${\bf c}=(c_{n})_{n\in \N}$ is a sequence of complex numbers. An important example is the M\"obius sequence. 
Recall that the M\"obius function $\mu (n)$ is, by definition, $\mu(n)=1$ if $n=1$; $\mu(n)=(-1)^{r}$ if $n=p_{1}\cdots p_{r}$ for $r$ distinct prime numbers $\{p_{i}\}_{i=1}^{r}$; $\mu (n)=0$ if $p^{2}|n$ for some prime number $p$. 
The M\"obius sequence ${\bf u}=(\mu (n))_{n\in \N}$ is the one generated by the M\"obius function. 
Following the idea of Sarnak~\cite{Sa1,Sa2},  we have the following definition.

\medskip
\begin{defn}[Disjointness]~\label{lindisj}
We say the sequence ${\bf c}=(c_{n})$ is {\em linearly disjoint} from the flow $T$ if
\begin{equation}\label{disjointness}
   \lim_{N\to \infty} \frac{1}{N}\sum_{n=1}^N c_n f(T^{n}x) =0
\end{equation}
for any $f\in C(X, \C)$ and any $x\in X$.
\end{defn}

Sarnak's conjecture (see~\cite{Sa1,Sa2}) says that the M\"obius sequence is linearly disjoint from all zero entropy flows. 
In~\cite{FanJiang}, we introduce a new concept called an oscillating sequence for the purpose of the study of this conjecture. 

\medskip
\begin{defn}[Oscillation]~\label{os}
The sequence ${\bf c}=(c_{n})_{n\in \N}$ is said to be an oscillating sequence if there is a constant $\lambda >1$ such that 
\begin{equation}~\label{cn}
K=\lim_{N\to \infty} \frac{1}{N} \sum_{n=1}^{N} |c_{n}|^{\lambda} <\infty
\end{equation} and if 
\begin{equation}~\label{oseq}
\lim_{N\to \infty} \frac{1}{N} \sum_{n=1}^{N} c_{n}e^{2\pi i n t}=0, \;\; \forall \; 0\leq t<1.
\end{equation}
\end{defn}

We proved in~\cite{FanJiang} that any oscillating sequence is linearly disjoint from all minimally mean attractable (MMA) and minimally 
mean-L-stable (MMLS) flows.   
In the same paper, 
we further proved that flows defined by all $p$-adic polynomials of integral coefficients,  all $p$-adic rational maps with good reduction, 
all automorphisms of the $2$-torus with zero topological entropy, all diagonalizable affine maps of the $2$-torus with zero topological entropy, 
all Feigenbaum maps, and all  orientation-preserving circle homeomorphisms are MMA and MMLS. 
Due to Davenport's theorem~\cite{Da}, the M\"obius sequence ${\bf u}$ is an oscillating sequence. 
Therefore, we confirmed Sarnak's conjecture for these flows which form a large class of zero topological entropy flows. Recently, Huang, Wang, and Zhang~\cite{HWZ} generalized a MLS flow to an ergodic  flow with discrete spectrum for invariant measures and proved that Sarnak's conjecture holds for these flows.

However, consider
$$
\T=\R^{2}/\Z^{2}
$$
the $2$-torus. In~\cite{FanJiang}, we also showed a counter-example as follows. Let
$$
T_{A, \alpha}=A{\bf x}^{t}+{\bf a}^{t}
$$ 
on the $2$-torus $\T$ where
$$
A= \left(
           \begin{array}{cc}
             1 & 0\\
             1 & 1 \\
           \end{array}
         \right)
$$
and
${\bf a}=(\alpha,0)\in \T$ is a non-zero constant point and ${\bf x} =(x,y)\in \T$ is a variable. Here ${\bf x}^{t}$ means the transpose of ${\bf x}$, that is,
$$
{\bf x}^{t}=
\left(
  \begin{array}{c}
    x \\
    y \\
  \end{array}\right).
$$
Then $T_{A, \alpha}$ is not MMLS on $\Td$ (it is MMA).
But Liu and Sarnak in~\cite{LS} and Wang in~\cite{W} showed that 
the M\"obius sequence ${\bf u}$ is linearly disjoint from this flows. Therefore, only the oscillation property is not enough for the purpose of the study of Sarnak's conjecture. We need the higher order oscillation property. There are two versions of a definition of the higher order oscillation (refer to~\cite[Remark 8]{FanJiang}).

\medskip
\begin{defn}[Weaker Version of Higher Order Oscillation]~\label{hoosw}
We call the sequence 
${\bf c}=(c_{n})_{n\in \N}$ 
a higher order oscillating sequence of order $m\geq 2$ if it satisfies (\ref{cn}) and if 
\begin{equation}~\label{hoosweq}
\lim_{N\to \infty}\frac{1}{N} \sum_{n=1}^{N} c_{n}e^{2\pi i n^{k}t}=0, \;\;\forall\; 1\leq k\leq m,\;\;\forall\; 0\leq t\leq 1.
\end{equation}
\end{defn}

Thanks to Hua's result~\cite{Hua}, we knew that the M\"obius sequence ${\bf u}$ is a higher order oscillating sequence of order $m$ for all $m\geq 2$ in this weaker version of the definition (Definition~\ref{hoosw}). Actually, according to~\cite{LZ}, we have that for any $A>0$,
\begin{equation}~\label{hoosd}
\sum_{n=1}^{N} c_{n}e^{2\pi i n^{k}t}=O_{A}\Big( N (\log N)^{-A}\Big), \;\;\forall\; 1\leq k\leq m,\;\;\forall\; 0\leq t\leq 1.
\end{equation}

\medskip
\begin{defn}[Stronger Version of Higher Order Oscillation]~\label{hooss}
We call the sequence 
$
{\bf c}=(c_{n})_{n\in \N}
$ 
a higher order oscillating sequence of order $m\geq 2$ if it satisfies (\ref{cn}) and if 
\begin{equation}~\label{hoosseq}
\lim_{N\to \infty} \frac{1}{N} \sum_{n=1}^{N} c_{n}e^{2\pi i P(n)}=0
\end{equation} 
for every real coefficient polynomial $P$ of degree $\leq m$. 
\end{defn}

In their paper~\cite[Lemma 2.1]{LS}, Liu and Sarnak showed that the M\"obius sequence ${\bf u}$
is also an higher order oscillating sequence of order $m$ for all $m\geq 2$ in this stronger version of the definitin (Definition~\ref{hooss}). They actually showed an estimation like the one in (\ref{hoosd}).

In a recent work~\cite{AJ}, we found another kind of higher order oscillating sequences of order $m$ for any $m\geq 2$ in the stronger version of the definition (Definition~\ref{hooss}), which is different from the one defined by an arithmetic function.

\medskip
\begin{thm}[Uniform Distribution~\cite{AJ}]~\label{unif}
Suppose $g$ is a positive $C^{2}$ function on $(1,\infty)$ with non-negative first and second derivatives. For a fixed real number $\alpha\neq 0$ and almost all real numbers $\beta>1$ (alternatively, for a fixed real number $\beta >1$ and almost all real number $\alpha$),
sequences 
$$
{\bf c} =\big( e^{2\pi i \alpha \beta^n g(\beta)}\big)_{n\in \N}
$$
are higher order oscillating sequence of order $m$ for any $m\geq 2$ in the stronger version of the definition (Definition~\ref{hooss}).
\end{thm} 

\begin{rem}~\label{more}
In particular, when $g\equiv 1$, sequences in Theorem~\ref{unif} are 
$$
{\bf c} =\big( e^{2\pi i \alpha \beta^n} \big)_{n\in \N}.
$$
\end{rem}

Since the weaker version of the definition (Definition~\ref{hoosw}) is more nature for the oscillation property in ergodic theory, 
therefore, it becomes an interesting question that can the weaker version of the definition (Definition~\ref{hoosw}) implies the stronger version of the definition (Definition~\ref{hooss})?
A proof of this will also provide another detailed proof of~\cite[Lemma 2.1]{LS}.

In this paper, 
we first give a definition of a distal flow on the $d$-torus (Definition~\ref{adf}). And then prove our main result in this paper that any higher order oscillating sequence of order $d$ is linearly disjoint 
from all affine distal flows on the $d$-torus. To present a clear idea, we first state and prove the main result for $d=2$ 
(Theorem~\ref{main1}) and then prove a consequence that
any higher order oscillating sequence of order $2$ is linearly disjoint 
from all affine flows on the $2$-torus with zero topological entropy (Corollary~\ref{cor1}). 
In particular, this reconfirms Sarnak's conjecture for all affine flows on the $2$-torus with zero topological entropy (Corollary~\ref{cormob2}). 
After that we state and prove 
the main result for $d>2$.
Thus it confirms Sarnak's conjecture 
for all affine distal flows on the $d$-torus  (Corollary~\ref{cormobd}).

\section{Statements of the Main Result}

Suppose $d\geq 2$ is an integer. Let 
$$
\Td=\R^{d}/\Z^{d}
$$
be the $d$-torus. Let $GL(d, \Z)$ be the space of all $d\times d-$matrices $A$ of integer entries with determinants $\det(A)=\pm 1$.
Let ${\bf x}=(x_{1}, \cdots, x_{d})\in \Td$ be a variable and denote
$$
{\bf x}^{t}=
\left(
  \begin{array}{c}
    x_{1} \\
    \vdots \\
    x_{d}
  \end{array}\right)
$$
the transpose of ${\bf x}$.
Then $A{\bf x}^{t}: \Td\to \Td$ is an automorphism of $\Td$. For ${\bf a}\in \Td$, we have an affine map $A{\bf x}^{t}+{\bf a}^{t}: \Td\to \Td$.
In order for this flow has entropy zero, the absolute values $\alpha$ of all eigenvalues of $A$ must be all $1$.
 If $\alpha$ is a complex number, then its complex conjugacy $\overline{\alpha}$ is also an eigenvalue.
 If $\alpha =e^{2\pi i \theta}$ for an irrational number $\theta$,  then the flow restricted on the union of the real part of the eigenspaces of $\alpha$ and $\overline{\alpha}$ is some kind rotation. Thus it is MMA and MMLS (see~\cite{FanJiang} or the proof of Corollary~\ref{cor1} when $d=2$).
Most interesting dynamics of this map is on the eigenspaces of $1$ and $-1$.   (If $\alpha =e^{2\pi i \theta}$ for a rational number $\theta=p/q$, then $A^{q}$ has $1$ or $-1$ as an eigenvalue.)  In this case,  except for the identity and the negative identity, by adding a non-zero shift ${\bf a}\in \Td$, we get a distal flow, which has certain polynomial expansion but still keeps entropy zero. However, for the notational simplicity, we include the identity and the negative identity in the following definition.

\medskip
\begin{defn}[Affine Distal Flow]~\label{adf}
We call an affine map
$$
T_{A,{\bf a}} =A{\bf x}^{t} +{\bf a}^{t}: \Td \to \Td
$$
an affine distal flow if all eigenvalues of $A$ are $1$ (or $-1$) and ${\bf a}\not ={\bf 0}\in \Td$.
\end{defn}

In the rest of the paper, we only use the stronger version of the definition (Definition~\ref{hooss}). 
In order to present our idea more clearly, we divide our main result in two cases: $d=2$ and $d>2$.
We state our main result in the case $d=2$ first.

\medskip
\begin{thm}[Main Theorem for $n=2$]~\label{main1}
Suppose ${\bf c}=(c_{n})_{n\in {\mathbb N}}$ is a higher order oscillating sequence of order $2$. Then it is linearly disjoint from all affine distal flows $T_{A,{\bf a}}$ on the $2$-torus $\T$.
\end{thm}

One of the consequences of this main result is that 

\medskip
\begin{cor}[Zero Entropy for $n=2$]~\label{cor1}
Suppose ${\bf c}=(c_{n})_{n\in {\mathbb N}}$ is a higher order oscillating sequence of order $2$. Then it is linearly disjoint from all affine flows $T_{A,{\bf a}}$ on the $2$-torus $\T$ with zero topological entropy.
\end{cor}

By combining Corollary~\ref{cor1} and~\cite[Lemma 2.1]{LS}, this reconfirms Sarnak's conjecture for all affine flows  with zero topological entropy.

\medskip
\begin{cor}[M\"obius Disjointness for $n=2$]~\label{cormob2}
The M\"obius sequence ${\bf u}=(\mu (n))_{n\in \N}$ is linearly disjoint from all affine flows $T_{A,{\bf a}}$ 
on the $2$-torus $\T$ with zero topological entropy.
\end{cor}

Now we state our main result in the case $d>2$.

\medskip
\begin{thm}[Main Theorem for Arbitrary $d>2$]~\label{main2}
For all $d> 2$, suppose ${\bf c}=(c_{n})_{n\in {\mathbb N}}$ is a higher order oscillating sequence of order $d$, 
then it is linearly disjoint from all affine distal flows $T_{A,{\bf a}}$ on the $d$-torus $\Td$.
\end{thm} 

By combining Theorem~\ref{main2} and~\cite[Lemma 2.1]{LS}, this confirms Sarnak's conjecture for 
all affine distal flows on the $d$-torus $\Td$.

\medskip
\begin{cor}[M\"obius Disjointness for $d>2$]~\label{cormobd}
For all $d> 2$, the M\"obius sequence ${\bf u}=(\mu (n))_{n\in \N}$ is linearly disjoint from all affine distal flows on the $d$-torus $\Td$.
\end{cor}

Our proof depends on the triangularization of an integral matrix. 
So we define a triangularizable affine distal flow. We say two affine distal flows $T_{A,{\bf a}}$ and $T_{B,{\bf b}}$ on the $d$-torus $\Td$ are topologically conjugate if there is a homeomorphism $h:\Td\to \Td$ such that
$$
T_{A,{\bf a}}\circ h=h\circ T_{B,{\bf b}}.
$$

\medskip
\begin{defn}~\label{tdfd}
We say a distal flow $T_{A,{\bf a}}$ is 
triangularizable if it is topologically conjugate to a distal flow 
$T_{B,{\bf b}}$
such that $B=(b_{ij})_{d\times d}$ is a upper-triangle matrix, that is, $b_{ij}=0$ for all $1\leq j<i\leq d$, and $b_{ii}=1$ for all $1\leq i\leq d$ (or $b_{ii}=-1$ for all $1\leq i\leq d$) and $b_{ij}\in \Z$ for all $1\leq i<j\leq d$.
\end{defn}

More precisely, the $d\times d-$matrix $B$ in Definition~\ref{tdfd} has the form
\begin{equation}~\label{form}
B= \pm \left(
\begin{array}{ccccccccccc}
1&b_{12}&b_{13}&\cdots&b_{1(d-1)}&b_{1d}\\
0&1&b_{21}&\cdots&b_{2(d-1)}&b_{2d}\\
\vdots&\vdots&\vdots&\vdots&\vdots&\vdots\\
0&0&0&\cdots&1&b_{(d-1)d}\\
0&0&0&\cdots&0&1
\end{array}\right)
\end{equation}
    
Before to prove Theorem~\ref{main2}, we will first prove the following lemma.

\begin{lem}~\label{tri}
For any $d\geq 2$, any distal flow $T_{A,{\bf a}}$ on the $d$-torus $\Td$ is triangularizable.
\end{lem}

\section{Proof of Theorem~\ref{main1}.}

Consider the space $C(\T, \C)$ of all complex continuous functions $f: \T \to \C$ with the supreme norm,
$$
\| f\|=\sup_{{\bf x}\in \T} |f({\bf x})|.
$$
Let ${\bf k}=(k,l)\in \Z^{2}$ and ${\bf x}=(x,y)$, define
$$
e({\bf k}\cdot {\bf x}) = e^{2\pi i (kx+ly)}, \quad {\bf x}\in \T.
$$
Then $e({\bf k}\cdot {\bf x})\in  C(\T, \C)$. From the Stone-Weierstrass theorem (refer to~\cite{R}), we have that 

\medskip
\begin{lem}~\label{sw}
The set 
$$
S=\Big\{ e({\bf k}\cdot {\bf x})\Big\}_{{\bf k}\in \Z^{2}}
$$ 
forms a dense subset in $C(\T, \C)$.
\end{lem}  

A linear combination $p$ of elements in $S$ is called a trigonometric polynomial. We can write $p$ as 
$$
p({\bf x}) = \sum_{k_{1}\leq k\leq k_{2}}\sum_{l_{1}\leq l\leq l_{2}} a_{kl} e^{2\pi i (kx +ly)}.
$$
Lemma~\ref{sw} implies that for any $f\in C(\T, \C)$, we have a sequence of trigonometric polynomials 
\begin{equation}~\label{stp2}
p_{q}({\bf x}) = \sum_{k_{1,q}\leq k\leq k_{2,q}} \sum_{l_{1,q}\leq l\leq l_{2,q}} a_{kl, q} e^{2\pi i (kx+ly)}.
\end{equation}
such that $\|f-p_{q}\|\to 0$ as $q\to \infty$. The sequence $\{p_{q}\}_{q\in\N}$ is called the trigonometric approximation of $f$.

Given $f\in C(\T, \C)$ and ${\bf x}\in \T$, let
$$
S_{N}({\bf x}) =  \frac{1}{N} \sum_{n=1}^{N} c_{n} f(T_{A, {\bf a}}^{n} {\bf x}).
$$
Let $\{ p_{q}\}_{q\in \N}$ be a sequence of trigonometric polynomials approximating $f$ in the supremum norm on $C(\T,\C)$. 

For any $\epsilon>0$, we have an integer $r>0$ such that
$$
\| f-p_{r}\| <\frac{\epsilon}{2K^{\frac{1}{\lambda}}}.
$$
Then 
$$
|S_{N}({\bf x})| \leq \Big| \frac{1}{N} \sum_{n=1}^{N} c_{n} | f(T_{A, {\bf a}}^{n} {\bf x}) -p_{r} (T_{A, {\bf a}}^{n} {\bf x})\Big| +\Big|\frac{1}{N} \sum_{n=1}^{N} c_{n} p_{r} (T_{A, {\bf a}}^{n} {\bf x})\Big| =I+II.
$$

For the estimation of $I$, we apply the H\"older inequality, 
$$
I\leq \Big( \frac{1}{N} \sum_{l=1}^{N} |c_{n}|^{\lambda}\Big)^{\frac{1}{\lambda}} 
 \Big( \frac{1}{N} \sum_{n=1}^{N} |f(T_{A, {\bf a}}^{n} {\bf x})-q_{r}(T_{A, {\bf a}}^{n} {\bf x})|^{\lambda'}\Big)^{\frac{1}{\lambda'}},
 $$
 where $\lambda'>1$ is the dual number of $\lambda$, that is, $1/\lambda +1/\lambda'=1$. 
 Thus we have 
$$
I\leq  K^{\frac{1}{\lambda}} \cdot \frac{\epsilon}{2K^{\frac{1}{\lambda}}} =\frac{\epsilon}{2}.
$$

For the estimation of $II$,  we first prove that $T_{A, {\bf a}}$ is triangularizable. Here we give a simple proof only working for  $d=2$ by using complex analysis. We give a proof for the general case in Lemma~\ref{tri}.  However, the proof we will give below must be used interestingly in the last step in the proof of Lemma~\ref{tri}. 

Suppose 
$$
  A = \left(
           \begin{array}{cc}
             a & b \\
             c & d \\
           \end{array}
         \right) \in GL(2,\Z)
$$
and all eigenvalues of $A$ are $1$. (If all eigenvalues are $-1$, then we consider $-A$.)  
It corresponds to the M\"obius transformation
$$
M(z) =\frac{az+b}{cz+d}, \quad z\in \C.
$$

Since $trace(M)=2$, it is a parabolic M\"obius transformation. So it has only one fixed point
$$
\frac{a-d}{2c}=\frac{p}{q}, \quad (p,q)=1,
$$
which is a rational point in $\C$. 
We have two integers $r$ and $s$ such that $pr-qs=1$ because of  the B\'ezout theorem.

Let 
$$
N(z) =\frac{pz+s}{qz+r}.
$$
It corresponds to the invertible integeral matrix
$$
  P= \left(
           \begin{array}{cc}
             p & s\\
             q & r\\
           \end{array}
         \right) \in GL(2,\Z)
$$
with $\det(P)=1$.
The M\"obius transformatiom $N$ 
maps $\infty$ to $p/q$ and has the inverse
$$
N^{-1}(z) =\frac{rz-s}{-qz+p}
$$
corresponding to the matrix
$$
  P^{-1}= \left(
           \begin{array}{cc}
             r & -s\\
             -q & p\\
           \end{array}
         \right) \in GL(2,\Z)
$$
with $\det(P^{-1})=1$.
Now consider the M\"obius transformation $N^{-1}\circ M\circ N$.
It is still a parabolic one and all coefficients are integers.
Most important, it only fixes $\infty$ and has no other fixed point in $\C$. Therefore
$$
 N^{-1}\circ M\circ N (z) = z+ t, \quad t\in \Z,
$$
which corresponds to the matrix 
\begin{equation}~\label{nf}
     B= \left(
           \begin{array}{cc}
             1 & t \\
             0 & 1 \\
           \end{array}
         \right), \quad t\in \Z.
\end{equation}
Thus, we have that
$$
P^{-1}AP= B.
$$
The map 
$$
h({\bf x}) =P{\bf x}^{t}:\T\to \T
$$
is a homeomorphism of the $2$-torus $\T$ and conjugates $T_{A,{\bf a}}$ to $T_{B, {\bf b}}$ for ${\bf b}=P^{-1} {\bf a} \not ={\bf 0}\in \T$, that is,
$$
T_{A, {\bf a}} \circ h = h\circ T_{B, {\bf b}}.
$$

Now let us continue to estimate $II$ under the assumption that $A$ is of the form (\ref{nf}).
Suppose 
${\bf a}=(a,b)$.
For any ${\bf x}=(x,y)\in \T$, let
$$
{\bf x}^{t}_{n}=T^{n}_{A,{\bf a}}{\bf x}^{t}=(x_{n}, y_{n}).
$$
Due to the fact that $A$ is of form (\ref{nf}), we have that
$$
y_{n}=y+nb
$$
and 
$$
x_{n} =x+(ty+a)n+ tb(\sum_{j=1}^{n-1} j)=x+(ty+a) n +\frac{tb}{2} n(n-1).
$$

Now using the formula (\ref{stp2}) for $p_{r}$,
we have that
$$
p_{r} (T_{A, {\bf a}}^{n} {\bf x}) =  
\sum_{k_{1} \leq k\leq k_{2}}\sum_{l_{1} \leq l\leq l_{2}} a_{kl,r} e^{2\pi i P_{kl,r}(n)}
$$
where
$$
P_{kl,r} (n) = \frac{tbk}{2} n^{2} +\Big( k \big( ty+a-\frac{tb}{2}\big) +lb\Big) n +(kx +ly)
$$
is a real coefficient polynomial of degree $2$ or $1$ or $0$.

Thus for the estimation of $II$, we have
$$
II =\Big|\frac{1}{N} \sum_{n=1}^{N} c_{n} \sum_{k_{1} \leq k\leq k_{2}}\sum_{l_{1}\leq l\leq l_{2}}  a_{kl,r} e^{2\pi i P_{kl, r} (n)}\Big| 
$$
$$
= \Big |\sum_{k_{1}\leq k\leq k_{2}}\sum_{l_{1}\leq l\leq l_{2}}  a_{kl,r} \frac{1}{N} \sum_{n=1}^{N}  c_{n} e^{2\pi i P_{kl, r} (n)} \Big|.
$$
Let 
$$
L=\max \{ |k_{1}|, |k_{2}|, |l_{1}|, l_{2}|, |a_{kl,r}| \; |\;  k_{1} \leq k\leq k_{2}, l_{1}\leq l\leq l_{2}\}.
$$
Since ${\bf c}=(c_{n})_{n\in \N}$ is a higher order oscillating sequence of order $2$, we can find an integer $M>r$ such that  for $N>M$,
$$
\Big|\frac{1}{N} \sum_{n=1}^{N}  c_{n} e^{2\pi i P_{kl,r} (n)} \Big|< \frac{\epsilon}{2L^{3}}, \quad \forall k_{1} \leq k\leq k_{2},\; l_{1}\leq l\leq l_{2}.
$$
This implies that 
$$
II< \frac{\epsilon}{2}.
$$
Therefore, we get that for all $N>M$,
$$
|S_{N}({\bf x})|< \epsilon.
$$
This says that $\lim_{N\to \infty} S_{N}({\bf x}) =0$. 
We proved Theorem~\ref{main1}.

\section{Proof of Corollary~\ref{cor1}.}

Let $\alpha$ and $\overline{\alpha}$ be two eigenvalues of $A$ in the complex field $\C$ and suppose $|\alpha|\geq 1$. 
The topological entropy of $T_{A, {\bf a}}$ is then $h(T_{A, {\bf a}}) =\log |\alpha|$. So $h(T_{A, {\bf a}})=0$ is equivalent to say that $|\alpha|=1$.
If $\lambda=e^{2\pi i \theta}$ for some $0< \theta <1$ but $\theta \not= 1/2$ (or when $\alpha =1$ and the other eigenvalue is $-1$), then $A$ is diagonalizable in the complex field $\C$. 
As we have proved~\cite[Proposition 8]{FanJiang},  $T_{A, {\bf a}}$ is an equicontinuous flow. Thus ${\bf c}= (c_{n})_{n\in \N}$ is linearly disjoint from 
$T_{A, {\bf a}}$  following our result in~\cite[Corollary 2]{FanJiang} since a higher order oscillating sequence of order $2$ is also an oscillating sequence.
When all eigenvalues of $A$ are $1$ (or $-1$), then $T_{A, {\bf a}}$ is a distal flow. It is a consequence of Theorem~\ref{main1}. This completes the proof.

\section{Proof of Lemma~\ref{tri}.}

In the proof of Theorem~\ref{main1}, we already saw a complex analysis proof of that any distal flow on the $2$-torus is triangularizbale. But this proof only works for $d=2$, although it is simple and neat.  Here we give a proof for the general case.
However, the argument we gave in the proof of Theorem~\ref{main1} has to be used in the last step of this proof.  

Suppose $d\geq 2$.  Suppose all eigenvalues of $A$ are $1$. (If all eigenvalues are $-1$, then we consider $-A$.)   
Since $A$ is an integral matrix and $1$ is its only eigenvalue, $A{\bf x}^{t}={\bf x}^{t}$ has an integer solution. Suppose ${\bf v}=(v_{1}, \cdots, v_{d})\in \Z^{d}$ is a solution such that 
$$
gcd (v_{1}, \cdots, v_{d})=1.
$$
Then, there are at least two of $v_{i}$, $1\leq i\leq d$, which are relatively prime. Without loss of generality, we assume $(v_{1},v_{2})=1$. From the B\'ezout theorem, there are two integers $r$ and $s$ such that $v_{1}r-v_{2}s=1$.
Consider the matrix
\begin{equation}~\label{form1}
P_{1}= \left(
\begin{array}{cccccccc}
v_{1}&s&0&0&\cdots&0&0\\
v_{2}&r&0&0&\cdots&0&0\\
v_{3}&0&1&0&\cdots&0&0\\
\vdots&\vdots&\vdots&\vdots&\vdots&\vdots&\vdots&\\
v_{d}&0&0&0&\cdots&0&1
\end{array}\right)
\end{equation}
Then $P_{1}\in GL(d, \Z)$ with $\det(P_{1})=1$. 

Let ${\bf e}_{1}=(1, 0,\cdots, 0)\in \Z^{d}$. We have that
$$
A{\bf v}^{t}={\bf v}^{t}\quad \hbox{and}\quad P_{1}{\bf e}_{1}^{t}={\bf v}^{t}.
$$
Thus, 
$$
P^{-1}_{1}{\bf v}^{t}={\bf e}_{1}^{t}.
$$
All these imply that
\begin{equation}~\label{form2}
P_{1}^{-1}AP_{1} =
\left(
\begin{array}{cc}
1&{\bf b}_{1(d-1)}\\
{\bf 0}^{t}_{d-1}&A_{1}
\end{array}\right)
\end{equation}
where ${\bf 0}_{d-1}=(0,\cdots,0)\in \Z^{d-1}$ and ${\bf b}_{1(d-1)}=(b_{12}, \cdots, b_{1d})\in \Z^{d-1}$ and $A_{1}\in GL(d-1, \Z)$ with $\det(A_{1})=1$. All eigenvalues of $A_{1}$ are $1$. 

Repeat the above argument for $A_{1}$, we have a $\widetilde{P}_{1}\in GL(d-1,\Z)$ with $\det( \widetilde{P}_{1})=1$ such that 
$$
\widetilde{P}_{1}^{-1}A_{1}\widetilde{P}_{1}
=
\left(
\begin{array}{ccc}
1&{\bf b}_{2(d-2)}\\
{\bf 0}^{t}_{d-2}& A_{2}
\end{array}\right)
$$
where ${\bf 0}_{d-2}=(0,\cdots,0)\in \Z^{d-2}$ and ${\bf b}_{2(d-2)}=(b_{23}, \cdots, b_{3d})\in \Z^{d-2}$ and $A_{2}\in GL(d-2, \Z)$ with $\det(A_{2})=1$.
Let 
$$
P_{2}=
\left(
\begin{array}{cc}
1&{\bf 0}_{d-1}\\
{\bf 0}^{t}_{d-1}& \widetilde{P}_{1}
\end{array}\right).
$$
Then we have 
$$
(P_{1}P_{2})^{-1}AP_{1}P_{2}=
\left(
\begin{array}{cc}
1&{\bf b}_{1(d-1)}\\
{\bf 0}^{b}_{d-1}& \widetilde{P}_{1}^{-1}A_{1}\widetilde{P}_{1}
\end{array}\right).
$$
 
Inductively, we obtain $P_{1}$, $P_{2}$, $\cdots$, $P_{d-2}$ such that
$$
\widehat{P}= P_{1}P_{2}\cdots P_{d-2}\in GL(d, \Z)
$$
with $\det(\widehat{P})=1$ and such that 
\begin{equation}~\label{form6}
\widehat{P}^{-1}A\widehat{P}
 =\left(
\begin{array}{ccccccccccc}
1&b_{12}&\cdots&b_{1(d-2)}&b_{1(d-1)}&b_{1d}\\
0&1&\cdots&b_{2(d-2)}&b_{2(d-1)}&b_{2d}\\
\vdots&\vdots&\vdots&\vdots&\vdots&\vdots\\
0&0&\cdots&1&b_{(d-2)(d-1)}&b_{(d-2)d}\\
0&0&\cdots&0&a&b\\
0&0&\cdots&0&c&d\\\end{array}\right).
\end{equation}

Let 
$$
A_{d-2}
=
\left(
\begin{array}{cc}
a&b\\
c& d
\end{array}\right).
$$
It is a $2\times 2$-matrix in $GL(2, \Z)$ with $det(A)=1$. All eigenvalues of $A_{d-2}$ are $1$. Now we apply the argument in the proof of Theorem~\ref{main1} to get $\widetilde{P}_{d-1}\in GL(2,\Z)$ with $\det(\widetilde{P}_{d-1})=1$ such that 
$$
\widetilde{P}_{d-1}^{-1}A\widetilde{P}_{d-1}=
\left(
\begin{array}{cc}
1&b_{(d-1)d}\\
0& 1
\end{array}\right), \quad b_{(d-1)d}\in \Z.
$$
Define 
$$
P_{d-1} =
\left(
\begin{array}{cc}
I_{d-2}&0\\
0& \widetilde{P}_{d-1}
\end{array}\right),
$$
where $I_{d-2}$ is the $(d-2)\times (d-2)$ identity matrix. And define 
$$
P =\widehat{P}P_{d-1}= P_{1}P_{2}\cdots P_{d-2}P_{d-1}.
$$
We finally get that $B=P^{-1}AP$ is of the form (\ref{form}) with $+$.
Let ${\bf b} =P^{-1}{\bf a}\in \Td$. 
Define $h({\bf x}) =P{\bf x}^{t}: \Td\to \Td$. It is a homeomorphism of the $d$-torus $\Td$ and we have that
$$
T_{A, {\bf a}} \circ h=h\circ T_{B, {\bf b}}.
$$
We completed the proof.

\section{Proof of Theorem~\ref{main2}.}

After the proof of Lemma~\ref{tri}, most of the proof of Theorem~\ref{main2} is similar to that of Theorem~\ref{main1} except for the notation and the estimation of $II$.
However for the independence of two theorems, we run a full proof again.

Consider the space $C(\Td, \C)$ of all complex continuous functions $f: \Td \to \C$ with the supreme norm,
$$
\| f\|=\sup_{{\bf x}\in \Td} |f({\bf x})|.
$$
Let ${\bf k}=(k_{1}, \cdots,k_{d})\in \Z^{d}$, define
$$
e({\bf k}\cdot {\bf x}) = e^{2\pi i (k_{1}x_{1}+\cdots + k_{d}x_{d})}, \quad {\bf x}\in \Td.
$$
Then $e({\bf k}\cdot {\bf x})\in  C(\Td, \C)$. From the Stone-Weierstrass theorem (refer to~\cite{R}), we have that 

\medskip
\begin{lem}~\label{swd}
The set 
$$
S=\Big\{ e({\bf k}\cdot {\bf x})\Big\}_{{\bf k}\in \Z^{d}}
$$ 
forms a dense subset in $C(\Td, \C)$.
\end{lem}  

A linear combination $p$ of elements in $S$ is called a trigonometric polynomial. We can write $p$ as 
$$
p({\bf x}) = \sum_{m_{1}\leq k_{1}\leq s_{1}}\cdots \sum_{m_{d}\leq k_{d}\leq s_{d}} a_{\bf k} e^{2\pi i (k_{1}x_{1}+\cdots +k_{d}x_{d})}.
$$
Lemma~\ref{swd} implies that for any $f\in C(\Td, \C)$, we have a sequence of trigonometric polynomials 
\begin{equation}~\label{stp}
p_{q}({\bf x}) = \sum_{m_{1q}\leq k_{1}\leq s_{1q}}\cdots \sum_{m_{dq}\leq k_{d}\leq s_{dq}} a_{{\bf k}, q} e^{2\pi i (k_{1}x_{1}+\cdots +k_{d}x_{d})}.
\end{equation}
such that $\|f-p_{q}\|\to 0$ as $q\to \infty$. The sequence $\{p_{q}\}_{q\in\N}$ is called the trigonometric approximation of $f$.

Given $f\in C(\Td, \C)$ and ${\bf x}\in \Td$, let
$$
S_{N}({\bf x}) =  \frac{1}{N} \sum_{n=1}^{N} c_{n} f(T_{A, {\bf a}}^{n} {\bf x}).
$$
Let $\{ p_{q}\}_{q\in \N}$ be a sequence trigonometric polynomials approximating $f$ in the supremum norm on $C(\Td,\C)$. 

For any $\epsilon>0$, we have an integer $r>0$ such that
$$
\| f-p_{r}\| <\frac{\epsilon}{2K^{\frac{1}{\lambda}}}.
$$
Then 
$$
|S_{N}({\bf x})| \leq \Big| \frac{1}{N} \sum_{n=1}^{N} c_{n} | f(T_{A, {\bf a}}^{n} {\bf x}) -p_{r} (T_{A, {\bf a}}^{n} {\bf x})\Big| +\Big|\frac{1}{N} \sum_{n=1}^{N} c_{n} p_{r} (T_{A, {\bf a}}^{n} {\bf x})\Big| =I+II.
$$

For the estimation of $I$, we apply the H\"older inequality, 
$$
I\leq \Big( \frac{1}{N} \sum_{l=1}^{N} |c_{n}|^{\lambda}\Big)^{\frac{1}{\lambda}} 
 \Big( \frac{1}{N} \sum_{n=1}^{N} |f(T_{A, {\bf a}}^{n} {\bf x})-q_{r}(T_{A, {\bf a}}^{n} {\bf x})|^{\lambda'}\Big)^{\frac{1}{\lambda'}},
 $$
 where $\lambda'>1$ is the dual number of $\lambda$, that is, $1/\lambda +1/\lambda'=1$. 
 Thus we have 
$$
I\leq  K^{\frac{1}{\lambda}} \cdot \frac{\epsilon}{2K^{\frac{1}{\lambda}}} =\frac{\epsilon}{2}.
$$

For the estimation of $II$, we can assume $A$ is of the form in (\ref{form}) with $+$ due to Lemma~\ref{tri} .

 Suppose 
${\bf a}=(a_{1},\cdots, a_{d})$.
For any ${\bf x} =(x_{1},\cdots,x_{d})\in \Td$, let
$$
{\bf x}^{t}_{n}=T^{n}_{A,{\bf a}}{\bf x}^{t}.
$$
Denote
${\bf x}_{n}=(x_{1}^{n},\cdots, x_{d}^{n})\in\Td.$

To show a clear idea of our proof, we first assume $d=3$.
We have
$$
x_{3}^{n} = x_{3} + na_{3},
$$
$$
x_{2}^{n} = x_{2} +(b_{23} x_{3} +a_{2}) n+ b_{23}a_{3} \sum_{j=1}^{n-1} j,
$$
and
$$
x_{1}^{n} = x_{1} +(b_{12}x_{2}+a_{1})n+b_{12} (b_{23} x_{3} +a_{2}) \sum_{j=1}^{n-1} j +b_{13}b_{23}a_{3} \sum_{k=1}^{n-1}\sum_{j=1}^{k} j.
$$
So we see that $x_{3}^{n}$ is a polynomial of $n$ of degree at most $1$, $x_{2}^{n}$ is a polynomial of $n$ of degree at most $2$, and $x_{1}^{n}$ is a polynomial of $n$ of degree at most $3$.

In general, for $d>2$,
due to the fact that $A$ is of form (\ref{form}) with $+$, one can see that
$$
x_{d}^{n}=x_{d}+na_{d},
$$
$$
x_{d-1}^{n} =x_{d-1}+(b_{(n-1)n}x_{n} +a_{n-1}) n+b_{(n-1)n}a_{n}\sum_{j_{1}=1}^{n-1}j_{1}.
$$
In general
$$
x_{i}^{n} = x_{i}^{n-1} + b_{i(i+1)} x_{i+1}^{n} +\cdots+ b_{id}x_{d}^{n-1}.
$$
So for $x_{d-2}^{n}$, as in the case $d=3$, it contains some single sums and a double sum $\sum_{j_{2}=1}^{n-1}\sum_{j_{1}=1}^{j_{2}} j_{1}$, which is a polynomial of degree $3$.
More general, suppose, in $x_{d-j-1}^{n}$, there is a term containing a degree $j$ polynomial $p_{j} (n)$, then in $x_{d-j}$, there is a term containing 
$$
p_{j+1}(n)=\sum_{k=1}^{n-1}p_{j}(k).
$$
which is a degree $j+1$ polynomial.
Thus we have that
$$
{\bf x}_{n}=(P_{1d} (n), \cdots, P_{d1}(n)),
$$
where $P_{(d-j+1)j}(n)$ is a polynomial of degree at most $j$ 
for each $1\leq j\leq d$. Coefficients are all real and unchanged when $n>d$.

Now using the formula (\ref{stp}), we have that
$$
p_{r} (T_{A, {\bf a}}^{n} {\bf x}) =  
\sum_{m_{1r} \leq k_{1}\leq s_{1r}} \cdots \sum_{ m_{dr} \leq k_{d}\leq s_{dr}} a_{{\bf k},r} e^{2\pi i P_{{\bf k},r}(n)},
$$
where
$$
P_{{\bf k},r} (n) = c_{{\bf k},r, d} n^{d} +\cdots +c_{{\bf k},r, 1}  n +c_{{\bf k},r,0}
$$
is a real coefficient polynomial of degree at most $d$.

For the estimation of $II$, we have
$$
II =\Big|\frac{1}{N} \sum_{n=1}^{N} c_{n} \sum_{m_{1r} \leq k_{1}\leq s_{1r}} \cdots \sum_{m_{dr}\leq k_{d}\leq s_{dr}}  a_{{\bf k},r} e^{2\pi i P_{{\bf k}, r} (n)}\Big| 
$$
$$
= \Big |\sum_{m_{1r}\leq k_{1}\leq s_{1r}}\cdots \sum_{m_{dr}\leq k_{d}\leq s_{dr}}  a_{{\bf k},r} \frac{1}{N} \sum_{n=1}^{N}  c_{n} e^{2\pi i P_{{\bf k}, r} (n)} \Big|.
$$
Let 
$$
L=\max \{ |m_{1r}|, \cdots, |m_{dr}|, |s_{1r}|, \cdots, |s_{dr}|, |a_{{\bf k},r}| \; |\;  m_{jr} \leq k_{j}\leq s_{jr}, 1\leq j\leq d\}.
$$
Since ${\bf c}=(c_{n})_{n\in \N}$ is a higher order oscillating sequence of order $d$, we can find an integer $M>r$ such that  for $N>M$,
$$
\Big|\frac{1}{N} \sum_{n=1}^{N}  c_{n} e^{2\pi i P_{{\bf k},r} (n)} \Big|< \frac{\epsilon}{2L^{3}}, \quad \forall m_{1r} \leq k_{1}\leq s_{1r}, \cdots, m_{dr}\leq k_{n}\leq s_{dr}.
$$
This implies that 
$$
II< \frac{\epsilon}{2}.
$$
Therefore, we get that for all $N>M$,
$$
|S_{N}({\bf x})|< \epsilon.
$$
This says that $\lim_{N\to \infty} S_{N}({\bf x}) =0$. 
We proved Theorem~\ref{main2}.

\medskip
\medskip
\noindent {\em Acknowledgement.} This work was done when I visited the National Center for Theoretical Sciences (NCTS) 
at National Taiwan University during 2016. I would like to thank NCTS for its hospitality. I also like to thank Professors Jung-Chao Ban and Chih-Hung Chang and other audiences for their spending times patiently to listen my proofs in a series of lectures in NCTS.

\medskip
\medskip
\noindent Yunping Jiang: Department of Mathematics, 
Queens College of the City University of New York,
Flushing, NY 11367-1597 and 
Department of Mathematics
Graduate School of the City University of New York
365 Fifth Avenue, New York, NY 10016

\noindent Email:yunping.jiang@qc.cuny.edu

\end{document}